# Controlled diffusion processes[*,†]

**Vivek S. Borkar**

*School of Technology and Computer Science,
Tata Institute of Fundamental Research,
Homi Bhabha Rd., Mumbai 400005, India.
e-mail:* `borkar@tifr.res.in`

**Abstract:** This article gives an overview of the developments in controlled diffusion processes, emphasizing key results regarding existence of optimal controls and their characterization via dynamic programming for a variety of cost criteria and structural assumptions. Stochastic maximum principle and control under partial observations (equivalently, control of nonlinear filters) are also discussed. Several other related topics are briefly sketched.

## 1. Introduction

The research on controlled diffusion processes took root in the sixties as a natural sequel to the developments in deterministic optimal control on one hand and in Markov decision processes on the other. From the former it inherited the legacy of 'compactness - lower semi-continuity' arguments for existence of optima and the Hamilton-Jacobi and (Pontryagin) maximum principle approaches to sufficient, resp. necessary conditions for optimality. From the latter it inherited the basic problem formulations corresponding to different cost functionals and more importantly, the notions of adapted (more generally, non-anticipative) controls, noisy observations, etc., which are peculiar to the stochastic set-up. As the field matured, this union proved to be greater than the sum of its parts and has contributed not only to its parent disciplines, but also to the theory of nonlinear partial differential equations, mathematical finance, etc. In this survey I shall attempt to give a comprehensive, though not exhaustive overview of the main strands of research in controlled diffusion processes.

The survey is organized as follows: The next section sets up the basic framework and solution concepts, defines the different classes of admissible control processes, and lists the standard problems in stochastic control classified according to the cost functional. Section 3 describes some motivating examples. Section 4 surveys the key results concerning the existence of optimal policies under resp. complete and partial observations. Section 5 deals with the characterization of optimality. The latter sphere of activity is dominated by dynamic programming and this is reflected in my write-up as well - comparatively less space is devoted to the other important strand, viz., stochastic maximum principle, for which pointers to literature are provided for greater detail. Section 6 briefly describes the computational issues. Section 7 presents an assortment of special topics.

---

[*]This is an original survey paper.
[†]Research supported in part by a grant for 'Nonlinear Studies' from the Indian Space Research Organization and the Defense Research and Development Organization, Government of India, administered through the Indian Institute of Science.

Throughout the article, I have given a few representative references, making no effort to be exhaustive (which would be an impossible task anyway).

## 2. The Control Problems

### 2.1. Solution concepts

Throughout what follows, we denote by $\mathcal{P}(S)$ the Polish space of probability measures on a Polish space $S$ with Prohorov topology. $\mathcal{L}(Z)$ will correspondingly stand for '*the law of* (an $S$−valued random variable) $Z$', viewed as an element of $\mathcal{P}(S)$. Also, for any $f : \mathcal{R}^+ \to S$ and $I \subset \mathcal{R}^+$, $f(I)$ denotes the trajectory segment $\{f(t), t \in I\}$.

The basic object of study here will be the $d$−dimensional ($d \geq 1$) controlled diffusion process $X(\cdot) = [X_1(\cdot), \cdots, X_d(\cdot)]^T$ described by the stochastic differential equation

$$X(t) = X_0 + \int_0^t m(X(s), u(s))ds + \int_0^t \sigma(X(s), u(s))dW(s), \qquad (1)$$

for $t \geq 0$. Here:

1. for a compact metric 'control space' $U$, $m(\cdot, \cdot) = [m_1(\cdot, \cdot), \cdots, m_d(\cdot, \cdot)]^T : \mathcal{R}^d \times U \to \mathcal{R}^d$ is continuous and Lipschitz in the first argument uniformly with respect to the second,
2. $\sigma(\cdot, \cdot) = [[\sigma_{ij}(\cdot, \cdot)]]_{1 \leq i,j \leq d} : \mathcal{R}^d \times U \to \mathcal{R}^{d \times d}$ is Lipschitz in its first argument uniformly with respect to the second,
3. $X_0$ is an $\mathcal{R}^d$−valued random variable with a prescribed law $\pi_0$,
4. $W(\cdot) = [W_1(\cdot), \cdots, W_d(\cdot)]^T$ is a $d$−dimensional standard Brownian motion independent of $X_0$,
5. $u(\cdot) : \mathcal{R}^+ \to U$ is the 'control process' with measurable paths, satisfying the *non-anticipativity condition*: for $t > s \geq 0$, $W(t) - W(s)$ is independent of $\{X_0, W(y), u(y), y \leq s\}$. (In other words, $u(\cdot)$ does not anticipate the future increments of $W(\cdot)$.)

We shall say that (1) is *non-degenerate* if the least eigenvalue of $\sigma(\cdot, \cdot)\sigma^T(\cdot, \cdot)$ is uniformly bounded away from zero, *degenerate* otherwise. The two solution concepts for (1) that we shall consider are:

1. *Strong solution:* Here we assume $X_0, W(\cdot), u(\cdot)$ to be given on a prescribed probability space $(\Omega, \mathcal{F}, P)$ and consider the corresponding $X(\cdot)$ given by (1). That there will be an almost surely unique $X(\cdot)$ can be proved by standard arguments using the Ito-Picard iterations as in [93], Ch. 4.

2. *Weak solution:* Here we assume that only the law of the pair $(X(\cdot), u(\cdot))$ is prescribed and consider any $(X(\cdot), u(\cdot), W(\cdot), X_0)$ on some probability space conforming to the above prescription. 'Uniqueness' then is interpreted as uniqueness in law.

These are exact counterparts of the corresponding notions for uncontrolled diffusions. Define

$$Lf(x, u) \stackrel{def}{=} \langle \nabla f(x), m(x, u) \rangle + \frac{1}{2} \text{tr}\left(\sigma(x, u)\sigma^T(x, u)\nabla^2 f(x)\right) \quad (2)$$

for $f \in C^2(\mathcal{R}^d)$. We may write $L_u f(x)$ for $Lf(u, x)$, treating $u$ as a parameter. Let $\{\mathcal{F}_t\}$ denote the natural filtration of $(X(\cdot), u(\cdot))$, i.e., $\mathcal{F}_t \stackrel{def}{=}$ the completion of $\cap_{s>t}\sigma(X(y), u(y), y \leq s)$. The weak solution is then equivalent to the following 'martingale' formulation:

*For any bounded twice continuously differentiable $f : \mathcal{R}^d \to \mathcal{R}$ with bounded first and second order partial derivatives, $f(X(t)) - \int_0^t L_{u(s)} f(X(s)) ds, t \geq 0$, is a martingale w.r.t. $\{\mathcal{F}_t\}$.*

It helps to think of the strong solution as the engineer's world view wherein $W(\cdot)$ is the noise input to a black box along with the chosen input $u(\cdot)$, leading to the 'output' $X(\cdot)$. The weak solution on the other hand represents the statistician's viewpoint in which one 'fits' the equation (1) to the known processes $(X(\cdot), u(\cdot))$ with $W(\cdot)$ being the noisy 'residuals'.

## 2.2. Control classes

The class of $u(\cdot)$ enunciated above is the most general class of controls that we shall consider, to be referred to as *non-anticipative* controls. Let $\{\mathcal{F}_t^X\}$ denote the natural filtration of $X(\cdot)$. Obviously, $u(\cdot)$ is adapted to $\{\mathcal{F}_t\}$. We shall say that $u(\cdot)$ is a *feedback* control if it is also adapted to $\{\mathcal{F}_t^X\}$, i.e., $u(t)$ at each $t$ is a function of the observed trajectory $X([0,t])$. We shall say that it is a *Markov* control if in addition $u(t) = v(t, X(t)), t \geq 0$, for a measurable $v : \mathcal{R}^+ \times \mathcal{R}^d \to U$. Finally, we say that it is a *stationary Markov* control if $u(t) = v(X(t)), t \geq 0$, for a measurable $v : \mathcal{R}^d \to U$.

We shall also need the relaxation of the notion of control process $u(\cdot)$ above to that of a *relaxed* control process. Here we assume that $U = \mathcal{P}(U_0)$ where $U_0$ is compact metrizable (whence so is $U$) and $m_i(\cdot, \cdot), 1 \leq i \leq d$, are of the form

$$m_i(x, u) = \int \bar{m}_i(x, y) u(dy), \ 1 \leq i \leq d,$$

for some $\bar{m}_i : \mathcal{R}^d \times U_0 \to \mathcal{R}$ that are continuous and Lipschitz in the first argument uniformly w.r.t. the second. Similarly, $\sigma(\cdot, \cdot)$ will be assumed to be of the form $[[\sigma_{ij}(x,u)]] =$ the nonnegative definite square-root of

$$\int \bar{\sigma}(x,y)\bar{\sigma}(x,y)u(dy)$$

for $\bar{\sigma}: \mathcal{R}^d \times U_0 \to \mathcal{R}^{d \times d}$ satisfying continuity / Lipschitz conditions akin to those for $\sigma$. See [111], pp. 132-134, about a discussion of the choice of square-root in the uncontrolled case, similar remarks apply here. In addition, we assume that all functions of the form $f(x,u), x \in \mathcal{R}^d, u \in U$, appearing in the cost criteria described below are of the form

$$f(x,u) = \int \bar{f}(x,y)u(dy), \ 1 \leq i \leq d,$$

for some $\bar{f}: \mathcal{R}^d \times U_0 \to \mathcal{R}$ satisfying the same regularity or growth conditions as those stipulated for $f$. We may write $u(t) = u(t, dy)$ to underscore the fact that it is a measure-valued process. Then the original notion of $U_0$-valued control $u_0(\cdot)$ (say) corresponds to $u(t, dy) = \delta_{u_0(t)}(dy)$, the Dirac measure at $u_0(t)$, for all $t$. We call such controls as *precise* controls. Precise feedback, Markov or stationary Markov controls may be defined accordingly. Intuitively, relaxed control generalizes the notion of randomized controls in discrete time problems, but this interpretation has to be treated with care: unlike in the discrete time case, we cannot have independent randomization at each $t$, as that would lead to measurability problems. A better picture is to view $dt \times u(t, dy)$ as a measure on $\mathcal{R}^+ \times U_0$. The set of relaxed controls is then the closure under weak* topology of the measures $dt \times \delta_{u_0(t)}(dy)$ corresponding to precise controls. In this sense, relaxed controls achieve the compactification and convexification of precise controls, which in turn form a dense subset therein. Unless mentioned otherwise, we shall work with the relaxed control framework. This notion was introduced in deterministic control by L. C. Young [117] and generalized to the stochastic case by Fleming [44]. It is a genuine relaxation in the sense that the corresponding joint laws of $(X(\cdot), u(\cdot))$ contain those corresponding to precise controls as a dense subset and therefore for most cost functionals of interest, the infimum over the latter equals the infimum over the former. The latter is often a minimum thanks to the compactification implicit in the relaxation.

Given (1) with $u(\cdot)$ a relaxed control, one can replace $u(\cdot), W(\cdot)$ in it by a $\tilde{u}(\cdot), \tilde{W}(\cdot)$ where $\tilde{u}(\cdot)$ is feedback and $\tilde{W}(\cdot)$ is another standard Brownian motion. In fact, $\tilde{u}(\cdot)$ is defined simply by $\int f d\tilde{u}(t) = E[\int f du(t)|\mathcal{F}_t^X]$ for $f$ in a countable subset of $C_b(U_0)$ that separates points of $U = \mathcal{P}(U_0)$ and $t \geq 0$. Conversely, if $(X_0, W(\cdot))$ are given on a probability space $(\Omega, \mathcal{F}, P)$ and a weak solution $(X'(\cdot), u'(\cdot), W'(\cdot))$ of (1) is available on some probability space with $\mathcal{L}(X_0') = \pi_0$ and $u'(\cdot)$ feedback, then this can be replicated in law by an $(X(\cdot), u(\cdot), W(\cdot))$ on a possibly enlarged $(\Omega, \mathcal{F}, P)$ with $(X_0, W(\cdot))$ as given. See [21], p. 18, for details in case of $\sigma$ without explicit control dependence. Extension to the more general case discussed here is straightforward. While the more flexible notion of weak solutions is usually the preferred one in dealing with controlled diffusions, the foregoing allows us to go back and forth between the strong and weak formulations to some extent.

In the non-degenerate case, (1) has a unique strong solution for a Markov control $v$ when $\sigma(\cdot, v(\cdot))$ is Lipschitz [115], which in particular includes the case when there is no explicit control dependence in $\sigma$. The Lipschitz requirement on $\sigma(\cdot, v(\cdot))$ can be relaxed to mere measurability for one and two dimensional problems along the lines of [111], pp. 192-194. (These results have been established for the case of bounded coefficients, but can be extended to, say, a 'linear growth' condition using a standard localization argument.) Also, the resulting processes can be shown to be strong Feller. On the other hand, in the non-degenerate case (1) always has a unique weak solution for feedback controls when $\sigma$ does not have explicit control dependence and is Lipschitz [93]. If $\sigma$ does have explicit control dependence and the control is stationary Markov, existence of a solution can be established ([76], p. 86-91), but not its uniqueness [99]. See, however, the results of [67] which show that under the non-degeneracy hypothesis, the property of having a unique strong solution is generic in a precise sense. (See also [77] for some instances where uniqueness is available.) In the degenerate case, neither existence nor uniqueness of either weak or strong solution is assured for general measurable controls. Under continuity (resp. Lipschitz) condition on $m(\cdot, v(\cdot)), \sigma(\cdot, v(\cdot))$, existence (resp. existence and uniqueness) of weak (resp., strong) solutions can be established even in the degenerate case [111].

Much of the literature on controlled diffusions does not include control in the diffusion matrix $\sigma(\cdot)$. There are some nontrivial reasons for this. The first is that for stationary Markov controls $u(\cdot) = v(X(\cdot))$, one is in general obliged to consider at best measurable $v(\cdot)$. As mentioned above, for a merely measurable diffusion matrix, even in the non-degenerate case only the existence of a weak solution is available. If one further allows explicit time dependence in $\sigma$, either through the control or otherwise, Lebesgue continuity of transition probabilities can be a problem [43]. Also, for a relaxed control process $\mu \in U$ with $\tilde{\sigma}(x, \mu) \stackrel{def}{=} \int \sigma(x, \cdot) d\mu$, $L_\mu f$ above needs to be defined as

$$\langle \nabla f(x), m(x, \mu) \rangle + \frac{1}{2} \mathrm{tr} \left( \int \sigma(x, u) \sigma^T(x, u) d\mu(u) \nabla^2 f(x) \right)$$

and not as

$$\langle \nabla f(x), m(x, \mu) \rangle + \frac{1}{2} \mathrm{tr} \left( \tilde{\sigma}(x, \mu) \tilde{\sigma}^T(x, \mu) \nabla^2 f(x) \right),$$

which can lead to problems of interpretation. (In situations where one can show that an optimum precise control exists, one can work around this problem.) This is not to say that the case of control-dependent diffusion matrix has not been studied. There are several situations, such as mathematical finance, where the control dependence of $\sigma$ cannot simply be wished away. Hence there has been a large body of work on problems with control-dependent drift. For example, the p.d.e. issues related to the HJB equations we mention later for problems with control-dependent diffusion matrix have been studied extensively in [7], [76]. More recently, Chinese mathematicians working in this area have developed an impressive body of work for this general case [116].

### 2.3. Cost structures

Let $k, c \in C(\mathcal{R}^d \times U), h \in C(\mathcal{R}^d), g \in C(\mathcal{R}^d \times \mathcal{R}^d), q \in C(U \times U)$, be prescribed functions with at most linear growth in the space (i.e., $x \in \mathcal{R}^d$) variable. Also, $c \geq 0$. Furthermore, in continuation of our relaxed control framework, $k, c$ are of the form $k(x, u) = \int \bar{k}(x, y) u(dy), c(x, u) = \int \bar{c}(x, y) u(dy)$, resp., for suitable $\bar{k}, \bar{c} \in C(\mathcal{R}^d \times U_0)$. Some standard cost functionals are:

1. *Finite horizon cost:* For $T > 0$, minimize

$$E[\int_0^T e^{-\int_0^t c(X(s),u(s))ds} k(X(t), u(t))dt + e^{-\int_0^T c(X(s),u(s))ds} h(X(T))]. \tag{3}$$

   Here $c$ is the discount function (*discount factor* if it is constant), $k$ the so called 'running cost' function and $h$ the terminal cost function.

2. *Cost up to exit time:* For an open set $D \subset \mathcal{R}^d$ with a smooth boundary $\partial D$ (more generally, boundary satisfying the 'exterior cone condition') and $\tau \stackrel{def}{=} min\{t \geq 0 : X(t) \notin D\}$, minimize

$$E[\int_0^\tau e^{-\int_0^t c(X(s),u(s))ds} k(X(t), u(t))dt + e^{-\int_0^\tau c(X(s),u(s))ds} h(X(\tau))]. \tag{4}$$

3. *Infinite horizon discounted cost:* For $c(\cdot, \cdot) \geq \delta > 0$, minimize

$$E[\int_0^\infty e^{-\int_0^t c(X(s),u(s))ds} k(X(t), u(t))dt]. \tag{5}$$

   This is popular in business applications where discounting is a real phenomenon and not merely a mathematical convenience.

4. *Average or 'ergodic' cost:* Minimize

$$\limsup_{T \to \infty} \frac{1}{T} \int_0^T E[k(X(t), u(t))]dt \tag{6}$$

   (the average version), or a.s. minimize

$$\limsup_{T \to \infty} \frac{1}{T} \int_0^T k(X(t), u(t))dt \tag{7}$$

   (the 'almost sure' version). These are popular in engineering applications where transients are fast, hence negligible, and one is choosing essentially from among the attainable 'steady states'.

5. *Risk-sensitive cost:* Minimize

$$E[e^{\int_0^T k(X(t),u(t))dt + h(X(T))}], \tag{8}$$

   or

$$\limsup_{T \to \infty} \frac{1}{T} \log E[e^{\alpha \int_0^T k(X(t),u(t))dt}], \tag{9}$$

where $\alpha > 0$ is a parameter. This cost functional has the advantage of involving 'all moments' of the cost, which matters when mere mean can be misleading. It also arises naturally in finance applications where compounding effects inherent in the formulation lead to the exponentiation in the cost [17]. Risk-sensitive control also has interesting connections with 'robust' control theory [38].

6. *Controlled optimal stopping:* Minimize

$$E[\int_0^\tau e^{-\int_0^t c(X(s),u(s))ds} k(X(t),u(t))dt + e^{-\int_0^\tau c(X(s),u(s))ds} h(X(\tau))] \quad (10)$$

over both admissible $u(\cdot)$ *and* all stopping times $\tau \geq 0$. The 'finite horizon' variation of this replaces $\tau$ above by $\tau \wedge T$ for a given $T > 0$.

7. *Impulse control:* Here one is allowed to reset the trajectory at stopping times $\{\tau_i\}$ from $X(\tau_i-)$ (the value immediately before $\tau_i$) to a new (non-anticipative) value $X(\tau_i)$, resp., with an associated cost $g(X(\tau_i), X(\tau_i-))$. The aim is to minimize

$$E[\int_0^T e^{-\int_0^t c(X(s),u(s))ds} k(X(t),u(t))dt +$$
$$\sum_{\tau_i \leq T} e^{-\int_0^{\tau_i} c(X(s),u(s))ds} g(X(\tau_i), X(\tau_i-)) + e^{-\int_0^T c(X(s),u(s))ds} h(X(T))],$$
$$(11)$$

over admissible $u(\cdot)$, reset times $\{\tau_i\}$, and reset values $\{X(\tau_i)\}$. Assume $g \geq \delta$ for some $\delta > 0$ to avoid infinitely many jumps in a finite time interval.

8. *Optimal switching:* Here one is allowed to switch the control $u(\cdot)$ at stopping times $\{\tau_i\}$ from $u(\tau_i-)$ (the value immediately before $\tau_i$) to a new (non-anticipative) value $u(\tau_i)$, resp., with an associated cost $q(u(\tau_i), u(\tau_i-))$. The aim is to minimize

$$E[\int_0^T e^{-\int_0^t c(X(s),u(s))ds} k(X(t),u(t))dt +$$
$$\sum_{\tau_i \leq T} e^{-\int_0^{\tau_i} c(X(s),u(s))ds} q(u(\tau_i), u(\tau_i-)) + e^{-\int_0^T c(X(s),u(s))ds} h(X(T))],$$
$$(12)$$

over reset times $\{\tau_i\}$, and reset values $\{u(\tau_i)\}$. Assume $q \geq \delta$ for some $\delta > 0$ to avoid infinitely many switchings in a finite time interval.

Infinite horizon discounted or ergodic versions of impulse and switching control can also be considered (see, e.g., [96], [108]). The *hybrid control* problem studied in [29] combines the last two above and more.

## 3. Examples

Here we sketch in brief some recent applications of controlled diffusions from literature. The description is necessarily brief and the reader is referred to the original sources for more details.

1. *Forest harvesting problem [3]:*

   In this problem, the so called 'stochastic forest stand value growth' is described up to extinction time $\gamma$ by

   $$X(t) = x + \int_0^t \mu(X(s))ds + \int_0^t \sigma(X(s))dW(s) - \sum_{\tau_k \leq t} \zeta_k,$$

   where $\gamma = \inf\{t \geq 0 : X(t) \leq 0\}$ (possibly $\infty$) and the non-negative, non-anticipative random variables $\{\tau_k\}, \{\zeta_k\}$, are respectively the cutting times and the quantities cut at the respective cutting times. The aim is to maximize the forest revenue $E[\sum_{\tau_k < \gamma} e^{-r\tau_k}(X(\tau_k) - c)]$, where $c > 0$ is the reforestation cost and $r > 0$ the discount factor. This is an impulse control problem.

2. *Portfolio optimization [73]:*

   The wealth process in portfolio optimization satisfies the s.d.e.

   $$dX(t) = X(t)[(\pi(t)\mu(t) + (1 - \pi(t))r(t))dt + \pi(t)\sigma(t)dW(t),$$

   where $\mu(\cdot), \sigma(\cdot)$ are known and $\pi(\cdot)$ is the $[0,1]$−valued control process that specifies the fraction invested in the risky asset, the remaining wealth being invested in a bond. Here $r(\cdot)$ is a fluctuating interest rate process satisfying

   $$dr(t) = a(t)dt + bdW'(t).$$

   Both $a(\cdot), b$ are assumed to be known and $W'(\cdot)$ is a Brownian motion independent of $W(\cdot)$. The aim is to maximize $E[X(T)^\gamma]$ for some $T, \gamma > 0$. ([73] considers a somewhat more general situation.) An alternative 'mean-variance' formulation in the spirit of Markowitz seeks to maximize a linear combination of the mean and negative variance of $X(T)$ [119]. A 'risk-sensitive' version of the problem, on the other hand, seeks to maximize

   $$\liminf_{T \uparrow \infty} -\frac{2}{\theta T} \log E[e^{-(2/\theta)X(T)}].$$

   See [78] for a slightly more general formulation.

3. *Production planning [13]:*

   Consider a factory producing a single good. Let $y(\cdot)$ denote its inventory level as a function of time and $p(\cdot) \geq 0$ the production rate. $\xi$ will denote the constant demand rate and $y_1, p_1$ resp. the factory-optimal inventory level and production rate. The inventory process is modelled as the controlled diffusion

$$dy(t) = (p(t) - \xi)dt + \sigma dW(t),$$

where $\sigma$ is a constant. The aim is to minimize over non-anticipative $p(\cdot)$ the discounted cost

$$E[\int_0^\infty e^{-\alpha t}[c(p(t) - p_1)^2 + h(y(t) - y_1)^2]dt]$$

where $c, h$ are known coefficients for the production cost and the inventory holding cost, resp.

4. *Heavy traffic limits of queues [62]:*

The following control problem arises in the so called Halfin-Whitt limit of multi-type multi-server queues : Consider a system of $d$ customer classes being jointly served by $N$ identical servers, with $\lambda_i, \mu_i, \gamma_i$ denoting the respective arrival, service and per customer abandonment rates for class $i$. Let $z_i = (\lambda_i/\mu_i)/\sum_j (\lambda_j/\mu_j), 1 \leq i \leq d$. In a suitable scaled limit (the aforementioned Halfin-Whitt limit), the vector of total number of customers of various classes present in the system satisfies the controlled s.d.e.

$$dX(t) = b(X(t), u(t))dt + \Sigma dW(t),$$

where the $i$−th component of $b(x, u)$ is $b_i(x, u) = -\theta \mu_i - \gamma_i(x_i - u_i) - \mu_i u_i$ and $\Sigma = diag[\sqrt{2\mu_1 z_1}, \cdots, \sqrt{2\mu_d z_d}]$. The parameter $\theta$ has the interpretation as the excess capacity of the server pool in a suitable asymptotic sense. The action space is state-dependent and at $x$, is

$$U(x) = \{u \in \mathcal{R}^d : u \leq x, \sum_i u_i = (\sum_i x_i) \wedge 0\}.$$

The $i$−th component of the control, $u_i(t)$, will correspond to a scaled limit of the number of servers assigned to the class $i$ at time $t$. The aim is to minimize the cost

$$E[\int_0^\infty e^{-\alpha t} c(X(t), u(t))dt]$$

for a discount factor $\alpha > 0$, where $c(x, u) = \sum_i (h_i + \gamma_i p_i)(x_i - u_i)$. Here $h_i, p_i$ are resp. the holding cost and the abandonment penalty for class $i$.

## 4. Existence results

### 4.1. Complete observations

Early existence theory in controlled diffusions was clearly motivated by the existing 'compactness-continuity' arguments from deterministic optimal control. The latter were based on establishing the sequential compactness of attainable trajectories of the state-control pairs in an appropriate function space and then establishing the continuity (more generally, lower semi-continuity) of the cost functional on it, whence the minimum was guaranteed to be attained. The first extensions of this approach considered the non-degenerate case without explicit control dependence in $\sigma$, under complete observations (i.e., the process $X(\cdot)$ is observed by the controller) and the finite horizon cost. Thus the $\mathcal{L}(X(\cdot))$ restricted to a finite time interval were absolutely continuous w.r.t. the law of the corresponding zero drift process, with the Radon-Nikodym derivative given by the Girsanov theorem. Establishing uniform integrability of these Radon-Nikodym derivatives, one obtained their relative sequential compactness in the $\sigma(L_1, L_\infty)$ topology by the Dunford-Pettis theorem. After establishing that every limit point thereof in this topology was also a legal Girsanov functional for some controlled diffusion, this was improved to compactness [5], [6], [37]. (See [47], [80] for some precursors which use more restrictive hypotheses.) [18] gives an ingenious argument to improve this to the existence of optimal Markov controls. [81] took a different approach based on establishing compactness of laws of the controlled processes in the space of probability measures on the trajectory space. While this is completely equivalent to the above for the non-degenerate case with control-independent $\sigma$, it provided a more flexible technique insofar as it could be extended to the degenerate case, control-dependent $\sigma$, infinite dimensional problems, etc.

The existence of optimal Markov controls can be read off the above for the case $c(\cdot, \cdot) \equiv$ a constant, simply from the fact that the one dimensional marginals of any controlled diffusion can be mimicked by another with a Markov control. This was first proved for the non-degenerate case [19], [59] and later extended to the degenerate case [15]. See [98] for similar results. To handle more general costs, it helps to view them as expectations with respect to appropriately defined 'occupation measures'. For example, the infinite horizon discounted cost

$$E[\int_0^\infty e^{-\alpha t} k(X(t), u(t)) dt]$$

($\alpha > 0$) can be written as $\int \bar{k} d\mu$ where the 'discounted occupation measure' $\mu$ is defined by:

$$\int f d\mu \stackrel{def}{=} E[\int_0^\infty e^{-\alpha t} \int f(X(t), y) u(t, dy) dt]$$

for $f \in C_b(\mathcal{R}^d \times U_0)$. This, of course, depends on the initial law which is assumed to be fixed. The set of attainable $\mu$ can be shown to be convex compact and in the non-degenerate case, one can show that each element thereof can be realized by a stationary Markov control (i.e., each $\mu$ can be mimicked by a stationary Markov control). In view of the lower semi-continuity of the map $\mu \to \int \bar{k} d\mu$, the desired existence result follows. This approach was initiated in [19] and carried out further in [21]. (In fact, one can show that the extreme points of this set correspond to *precise* stationary Markov controls, see the discussion of the ergodic control problem below.) In the degenerate case, such a 'mimicry theorem' for occupation measures seems unavailable, but the existence of an optimal Markov (for finite horizon problems) or stationary Markov (for discounted infinite horizon problem or control up to exit time) controls can be established by adapting Krylov's Markov selection procedure ([111], Ch. 12). This was done in [40], [63] following a suggestion of Varadhan. Another variation of the above, applicable to the degenerate case, looks at equivalence classes of $\mathcal{L}(X(\cdot), u(\cdot))$ whose marginals agree a.e. and shows that the extremal equivalence classes in fact correspond to Markov controls [23]. See [24] for a further variation on this theme.

Throughout the foregoing, as one might expect, one has to weaken 'stationary Markov' to 'Markov' if the cost and / or the drift and / or the diffusion matrix of (1) have explicit time dependence. Also, for $U_0 \subset \mathcal{R}^m$, the compactness assumption on $U_0$ can be dropped by penalizing large $||u(t)||$, e.g., by including the term $\frac{1}{2}||u(t)||^2$ in the running cost.

The occupation measure approach is most successful for the ergodic control problem. This has been studied mostly for the case when $\sigma$ does not have explicit control dependence, because of the possible non-uniqueness of solutions under stationary Markov controls when it does. (More generally, one would have to work with 'the set of all solutions' for a stationary Markov control rather than *the* solution.) Consider the non-degenerate case first. Let $v(\cdot)$ be a stationary Markov control such that the corresponding $X(\cdot)$ is positive recurrent and therefore has a unique stationary distribution $\eta^v \in \mathcal{P}(\mathcal{R}^d)$. Define the corresponding ergodic occupation measure as $\mu^v(dx, dy) \stackrel{def}{=} \eta^v(dx) v(x, dy)$. One can show that the set $G$ of attainable $\mu^v$'s is closed convex with its extreme points corresponding to precise stationary Markov controls. We can say much more: define the empirical measures $\{\nu_t\}$ by:

$$\int f d\nu_t \stackrel{def}{=} \frac{1}{t} \int_0^t \int f(X(s), y) u(s, dy) ds, \ f \in C_b(\mathcal{R}^d \times U_0), t > 0.$$

Let $\bar{\mathcal{R}} = \mathcal{R}^d \cup \{\infty\} \stackrel{def}{=}$ the one point compactification of $\mathcal{R}^d$ and view $\nu_t$ as a random variable in $\mathcal{P}(\bar{\mathcal{R}} \times U)$ that assigns zero mass to $\{\infty\} \times U$. Then as $t \to \infty$,

$$\nu_t \ \to \ \{\nu : \nu(A) = a\nu'(A \cap (\{\infty\} \times U)) + (1-a)\nu''(A \cap (\mathcal{R}^d \times U)) \ \forall \ A$$
$$\text{Borel in } \bar{\mathcal{R}} \times U_0, \text{ with } a \in [0,1], \ \nu' \in \mathcal{P}(\{\infty\} \times U), \nu'' \in G\}$$

almost surely. This allows us to deduce the existence of an optimal precise stationary Markov control for the 'a.s.' version of the ergodic control problem in two cases: $(i)$ Under a suitable 'stability' condition (such as a convenient 'stochastic Liapunov condition') that ensures compactness of $G$ and a.s. tightness of $\{\nu_t\}$, or $(ii)$ a condition that penalizes escape of probability mass to infinity, such as the 'near-monotonicity condition':

$$\liminf_{||x|| \to \infty} \min_u k(x, u) > \beta,$$

where $\beta \stackrel{def}{=}$ the optimal cost [25]. The latter condition is often satisfied in practice. The 'average' version of the ergodic cost can be handled similarly.

As always, the degenerate case is much harder. Here one shows that, as in the non-degenerate case, $G$ is characterized as $\{\mu : \int Lfd\mu = 0\}$ for $f \in$ a sufficiently rich class of functions in $C^2(\mathcal{R}^d)$. That a $\mu \in G$ would satisfy $\int Lfd\mu = 0$ is easy to see for the stipulated $f$. The hard part is the reverse implication: One shows that there exists a stationary pair $(X(\cdot), u(\cdot))$ that has $\mu$ as its marginal. This extends an important result of [39] to the controlled case [109]. See [15], [79] for some extensions. This characterization helps establish $G$ as a closed convex set, leading to existence of optimal ergodic pairs $(X(\cdot), u(\cdot))$ under suitable (somewhat stronger) stability or near-monotonicity conditions [21]. This can be refined to an optimal stationary Markov control by means of a limiting argument using Krylov selection for the discounted cost as the discount factor approaches zero [16].

It should be mentioned that in the non-degenerate case, one often has classical solutions to the associated HJB equation as we see later and the existence of an optimal precise stationary Markov (or Markov in the finite horizon case) control can be read off the HJB theory. Thus direct existence results described above at best give some additional insight, except in some 'non-classical' situations like the constrained problems we encounter later. The 'occupation measure' viewpoint above is also the basis for the linear programming approach we discuss later. In the degenerate case, however, there is significant motivation to pursue these.

Finally, we note that such 'direct' existence results are also possible for more general problems involving impulsive and switching controls, etc. See, e.g., [29].

## 4.2. Partial observations

This corresponds to the situation where there is another $m$−dimensional 'observations' process $Y(\cdot)$ given by

$$Y(t) = \int_0^t b(X(s))ds + W'(t), \ t \geq 0,$$

where $b : \mathcal{R}^d \to \mathcal{R}^m$ is Lipschitz and $W'(\cdot)$ is an $m$−dimensional standard Brownian motion independent of $W(\cdot)$. $W'(\cdot)$ corresponds to (integrated) 'observation noise', as opposed to the 'signal noise' $W(\cdot)$. The situation when the two are not independent is called the 'correlated noise' case and has also been studied in literature. The objective is to optimize one of the above cost functionals over all control processes $u(\cdot)$ adapted to the natural filtration of $Y(\cdot)$, denoted $\{\mathcal{F}_t^Y\}$. We shall call these *strict sense admissible controls*, to contrast them with *wide sense admissible* controls to be defined later.

The correct 'state' (or 'sufficient statistics') process in this case turns out to be the regular conditional law $\pi_t$ of $X(t)$ given $\mathcal{G}_t \stackrel{def}{=}$ the right-continuous completion of $\sigma(Y(s), u(s), s \leq t)$ for $t \geq 0$. (For strict sense admissible $u(\cdot)$, this is the same as $\{\mathcal{F}_t^Y\}$.) This is a $\mathcal{P}(\mathcal{R}^d)$−valued process whose evolution is described by the Fujisaki-Kallianpur-Kunita equation of nonlinear filtering described as follows: Let $\nu(f) \stackrel{def}{=} \int f d\nu$ for any non-negative measure $\nu$ on $\mathcal{R}^d$ and $f \in C_b(\mathcal{R}^d)$. Then for $f \in C_0^2(\mathcal{R}^d) \stackrel{def}{=}$ the space of twice continuously differentiable $f : \mathcal{R}^d \to \mathcal{R}$ which vanish at infinity along with their first and second order partial derivatives, one has

$$\pi_t(f) = \pi_0(f) + \int_0^t \pi_s(L_{u(s)}f)ds + \int_0^t \langle \pi_s(bf) - \pi_s(b)\pi_s(f), d\hat{Y}(s)\rangle. \quad (13)$$

Here the so called 'innovations process'

$$\hat{Y}(t) \stackrel{def}{=} Y(t) - \int_0^t \pi_s(b)ds, \ t \geq 0,$$

is an $m-$dimensional standard Brownian motion independent of $(X_0, W(\cdot))$ and generating the same filtration as $Y(\cdot)$ [1].

Let $\bar{\mathcal{F}}_t \stackrel{def}{=}$ the right-continuous completion of

$$\sigma(X(s), Y(s), u(s), W(s), W'(s), s \leq t)$$

for $t \geq 0$. Let $(\Omega, \mathcal{F}, P)$ denote the underlying probability triple where $\mathcal{F} = \vee_t \bar{\mathcal{F}}_t$ without loss of generality. Define a new probability measure $P_0$ on $(\Omega, \mathcal{F})$ by:

$$\frac{dP|_{\bar{\mathcal{F}}_t}}{dP_0|_{\bar{\mathcal{F}}_t}} = \Lambda_t \stackrel{def}{=} e^{\int_0^t \langle b(X(s)), dY(s)\rangle - \frac{1}{2}\int_0^t ||b(Y(s))||^2 ds}, \ t \geq 0.$$

By Girsanov's theorem, under $P_0$, $Y(\cdot)$ itself is an $m-$dimensional standard Brownian motion independent of $(X_0, W(\cdot))$. Define the process of *unnormalized conditional laws* $\mu_t, t \geq 0$, taking values in $\mathcal{M}(\mathcal{R}^d)$, the space of non-negative measures on $\mathcal{R}^d$ with the weak* topology, as follows:

$$\mu_t(f) \stackrel{def}{=} E_0[f(X(t))\Lambda_t|\mathcal{G}_t]$$

for a countable collection of $f \in C_b(\mathcal{R}^d)$ that separates points of $\mathcal{M}(\mathcal{R}^d)$, $E_0[\ \cdot\ ]$ being the expectation under $P_0$. This evolves according to the *Duncan-Mortensen-Zakai* equation

$$\mu_t(f) = \mu_0(f) + \int_0^t \mu_s(L_{u(s)}f)ds + \int_0^t \langle \mu_s(bf), dY(s)\rangle \quad (14)$$

for $f \in C_0^2(\mathcal{R}^d)$. This has the advantages of linearity and the fact that viewed under $P_0$, $Y(\cdot)$ itself is the driving Brownian motion. $\mu_t, t \geq 0$ is interconvertible with $\pi_t, t \geq 0$, through:

$$\pi_t(f) = \frac{\mu_t(f)}{\mu_t(\mathbf{1})},$$

$$\mu_t(f) = \pi_t(f) e^{\int_0^t \langle \pi_s(b), dY(s)\rangle - \frac{1}{2}\int_0^t ||\pi_s(b)||^2 ds}.$$

Here **1** is the constant function identically equal to 1. Thus $\mu_t$ is an equivalent state variable. The first equality above justifies the adjective 'unnormalized'.

Yet another equivalent state variable is defined by the process $\varphi_t, t \geq 0$, given by

$$\varphi_t(f) \stackrel{def}{=} \mu_t(e^{-\langle b, Y_t \rangle} f), \ f \in C_b(\mathcal{R}^d).$$

Thus $\mu_t(f) = \varphi_t(e^{\langle b, Y_t \rangle} f)$. Suppose that $b(\cdot)$ is twice continuously differentiable. By an 'integration by parts' argument, $\{\varphi_t\}$ is seen to evolve according to

$$\varphi_t(f) = \varphi_0(f) + \int_0^t \varphi(\tilde{L}_{u(s),s} f) ds, \tag{15}$$

for $f \in C_0^2(\mathcal{R}^D)$. Here,

$$\begin{aligned}\tilde{L}_{u,s} f(x) \stackrel{def}{=}\ & L_u f(x) - \langle \nabla f(x), \sigma(x)\sigma^T(x) Db^T(x) Y(s) \rangle \\ & + (\frac{1}{2} \langle Y(s), Db(x)\sigma(x)\sigma^T(x) Db^T(x) Y(s) \rangle \\ & - \langle Y(s), Db(x) m(x,u) + \ell(x) \rangle - \frac{1}{2} \|b(x)\|^2) f,\end{aligned}$$

where $Db$ is the Jacobian matrix of $b$ and $\ell_i(x) \stackrel{def}{=} \frac{1}{2}\mathrm{tr}(\sigma^T(x) \nabla^2 b_i(x) \sigma(x))$ for $1 \leq i \leq m$. (15) is an *ordinary* parabolic p.d.e. (as opposed to the *stochastic* p.d.e.s (13) and (14)) with the sample path of $Y(\cdot)$ appearing as a random parameter. Hence this is called the pathwise filter. Standard p.d.e. theory ensures the well-posedness of the pathwise filter, from which that of the Fujisaki-Kallianpur-Kunita and Duncan-Mortensen-Zakai filters may be deduced using the conversion formulas [65].

Note that for $f \in C_b(\mathcal{R}^d \times U)$,

$$\begin{aligned}E[f(X(t), u(t))] &= E[\pi_t(f(\cdot, u(t)))] \\ &= E_0[\mu_t(f(\cdot, u(t)))] \\ &= E_0[\varphi_t(e^{\langle Y(t), b(\cdot) \rangle} f(\cdot, u(t)))].\end{aligned}$$

Thus, for example, the finite horizon cost

$$E[\int_0^T k(X(s), u(s)) ds + h(X(T))] \tag{16}$$

equals

$$E[\int_0^T \pi_t(k(\cdot, u(s))) ds + \pi_T(h)] \tag{17}$$

or,

$$E_0[\int_0^T \mu_t(k(\cdot, u(s))) ds + \mu_T(h)] \tag{18}$$

or,

$$E_0[\int_0^T \varphi_t(e^{\langle Y(s), b(\cdot) \rangle} k(\cdot, u(s))) ds + \varphi_T(e^{\langle Y(T), b(\cdot) \rangle} h)]. \tag{19}$$

Hence the control problem of minimizing (16) under partial observations can be viewed as the equivalent problem of controlling the $\mathcal{P}(\mathcal{R}^d)$−valued (resp., $\mathcal{M}(\mathcal{R}^d)$−valued) process $\{\pi_t\}$ (resp., $\{\mu_t\}, \{\varphi_t\}$) with cost (17) (resp., (18), (19)). These are called *separated* control problems because they 'separate', i.e., compartmentalize the two issues of state estimation and control. When $\pi_t$ or any of the equivalent state variables can be characterized by finite dimensional 'sufficient statistics', this can be reduced to a finite dimensional control problem. Such instances are rare, but include the important case of linear systems with linear observations (i.e., $m, b$ are linear and $\sigma$ a constant), Gaussian $\pi_0$, and a quadratic cost. Here $\pi_t$ is Gaussian and is completely characterized by its first two moments, corresponding to first two *conditional* moments of the state given observations.

For discounted cost (3) with discount function $c$, one replaces $\pi_t$ by $\tilde{\pi}_t$ in (13) and (16), where $\tilde{\pi}(f) \stackrel{def}{=} E[e^{-\int_0^t c(X(s),u(s))ds} f(X(t))|\mathcal{G}_t]$ for $f \in C_b(\mathcal{R}^d)$. Correspondingly, replace $L_u$ by $L_u - c(\cdot, u)$ in (13). Similar adjustments can be made for (14) and (15).

The existence of optimal strict sense admissible control for this problem remains an open issue. The best known result is the existence of an optimal *wide sense admissible* control [48]. Say $u(\cdot)$ is wide sense admissible if for each $t \geq 0, Y(t+\cdot) - Y(t)$ is independent of $(X_0, W(\cdot), \{u(s), Y(s), s \leq t\})$ under $P_0$. This clearly includes strict sense admissible controls and can be shown to be a valid relaxation thereof in the sense that the infimum of the cost over either set is identical. The proof technique for the existence claim is based on weak convergence arguments for measure-valued processes akin to the complete observations case and exploits the fact that wide sense admissibility is preserved under convergence of joint laws of the processes concerned. A refinement based on Krylov's Markov selection procedure leads to the existence of an optimal Markovian control (i.e., one that at each time depends on the current value of the measure-valued filter) for the separated control problem [41].

Similar developments are possible for (5). For (4), one replaces $\{\pi_t\}$ with a suitably modified measure-valued process that is supported on $D$ [20]. For ergodic control, the separated control problem can be formulated for $\{\pi_t, u(t), t \geq 0\}$ exactly as above and the existence theory is analogous to that for the degenerate diffusions under complete observations described above modulo additional technicalities [15], [16]. For risk-sensitive control, one needs a modification of the measure-valued process along the lines described above for taking care of the discount function $c$.

## 5. Characterization of Optimal Controls

### 5.1. HJB equation - the classical case

We begin with the dynamic programming principle, which is usually the preferred approach to characterization of optimality in controlled diffusions (see, e.g., [95]). To start with, consider the non-degenerate case. Consider for example the finite horizon cost. Define the 'value function'

$$V(x,t) = \inf E[\int_t^T e^{-\int_t^y c(X(s),u(s))ds} k(X(y),u(y))dy$$
$$+ e^{-\int_t^T c(X(s),u(s))ds} h(X(T))|X(t) = x],$$

where the infimum is over all admissible controls. Then by the standard dynamic programming heuristic, for $t' > t$,

$$V(x,t) = \inf E[\int_t^{t'} e^{-\int_t^y c(X(s),u(s))ds} k(X(y),u(y))dy$$
$$+ e^{-\int_t^{t'} c(X(s),u(s))ds} V(X(t'),t')|X(t) = x].$$

In words, if one is at point $x$ at time $t$, then the 'minimum cost to go' is the minimum of the sum of the cost paid over $[t, t']$ plus the minimum cost to go from time $t'$ onwards. Let $t' = t + \Delta$ for some small $\Delta > 0$. Then

$$V(x,t) \approx \inf E[k(X(t),u(t))\Delta +$$
$$e^{-c(X(t),u(t))\Delta} V(X(t+\Delta),t+\Delta)|X(t) = x].$$

Thus

$$\inf E[k(X(t),u(t))\Delta + e^{-c(X(t),u(t))\Delta} V(X(t+\Delta),$$
$$t+\Delta) - V(x,t)|X(t) = x]/\Delta \approx 0.$$

Assuming sufficient regularity of $V$, letting $\Delta \to 0$ formally leads to

$$\frac{\partial V}{\partial t} + \min_u (k(x,u) + \langle \nabla_x V(x,t), m(x,u) \rangle - c(x,u)V(x,t)$$
$$+ \frac{1}{2}\text{tr}\left(\sigma(x,u)\sigma^T(x,u)\nabla_x^2 V(x,t)\right)) = 0, \quad (20)$$

where $\nabla_x, \nabla_x^2$ denote the gradient and the Hessian in the $x$ variable. This is the HJB equation for the finite horizon control problem, with the boundary condition $V(x,T) = h(x) \; \forall \; x$.

The above is an instance of how the dynamic programming heuristic is used to guess the correct HJB equation. The equation is then analyzed by invoking the standard p.d.e. theory. For example, under appropriate boundedness and regularity conditions on $m, \sigma, k, h, c$, (20) has a unique solution in the Sobolev space $W_p^{2,1}(\mathcal{R}^d \times [0,T])$ for any $p$, $2 \leq p < \infty$. When $\sigma$ does not have explicit control dependence, this is a *quasi-linear* p.d.e. as opposed to a fully nonlinear one and the existence of a unique solution can be established in the class of bounded $f : \mathcal{R}^d \times \mathcal{R}^+ \to \mathcal{R}$ which are twice continuously differentiable in the first variable and once continuously differentiable in the second [86]. In either case, that this solution indeed equals the value function follows by a straightforward argument based on Krylov's extension of the Ito formula ([76],

p. 122). Implicit in this is the following 'verification theorem': A Markov control $v : \mathcal{R}^d \times [0,T] \to U$ is optimal if and only if

$$v(x,t) \in \mathrm{Argmin}_u \left( L_u V(x,t) + k(x,u) \right) \quad \text{a.e.} \tag{21}$$

The existence of a measurable $v(\cdot)$ satisfying the above follows from a standard measurable selection theorem [102].

For control up to exit time, it makes obvious sense to define

$$V(x) \stackrel{def}{=} \inf E\left[ \int_0^\tau e^{-\int_0^t c(X(s),u(s))ds} k(X(t),u(t)) dt \right.$$
$$\left. + e^{-\int_0^\tau c(X(s),u(s))ds} h(X(\tau)) | X(0) = x \right],$$

the infimum being over all admissible controls. There is no explicit time dependence in $V$ because the 'possible futures till $\tau$' look the same from a given state regardless of when one arrived there. A heuristic similar to the above leads to the HJB equation

$$\min_u (k(x,u) - c(x,u)V(x) + LV(x,u)) = 0 \tag{22}$$

with $V(x) = h(x)$ for $x \in \partial D$. This has a unique solution in $W^2_{p,loc}(D) \cap C(\bar{D})$ [53]. A verification theorem for optimal *stationary* Markov controls along the lines of (21) can be established.

For the infinite horizon discounted cost, the HJB equation for

$$V(x) \stackrel{def}{=} \inf E\left[ \int_0^\infty e^{-\int_0^t c(X(s),u(s))ds} k(X(t),u(t)) dt | X(0) = x \right] \tag{23}$$

is (22) on the whole space and for $k$ bounded from below, the 'value function' defined above is its least solution in $W^2_{p,loc}(\mathcal{R}^d)$. An appropriate verification theorem holds. In both this and the preceding case, '$W^2_{p,loc}$' can be replaced by '$C^2$' in the quasi-linear case corresponding to control-independent $\sigma(\cdot)$.

The situation for ergodic control is more difficult. Let $V^\alpha$ denote the $V$ of (23) when $c \equiv$ a constant $\alpha > 0$. Define $\bar{V}^\alpha \stackrel{def}{=} V^\alpha - V^\alpha(0)$. Then $\bar{V}^\alpha$ satisfies

$$\min_u (k(x,u) - \alpha \bar{V}^\alpha(x) - \alpha V^\alpha(0) + L\bar{V}^\alpha(x,u)) = 0. \tag{24}$$

Under suitable technical conditions (such as near-monotonicity or stability conditions mentioned above) one can show that as $\alpha \downarrow 0$, $\bar{V}^\alpha(\cdot)$ and $\alpha V^\alpha(0)$ converge along a subsequence to some $V, \beta$ in resp. an appropriate Sobolev space and $\mathcal{R}$. Letting $\alpha \downarrow 0$ along this subsequence in (24), these are seen to satisfy

$$\min_u (k(x,u) - \beta + LV(x,u)) = 0.$$

This is the HJB equation of ergodic control. One can show uniqueness of $\beta$ as being the optimal ergodic cost and of $V$ up to an additive scalar in an appropriate function class depending on the set of assumptions one is working with. A verification theorem holds [9], [25].

For risk-sensitive control, the HJB equations are

$$\min_u(\frac{\partial V}{\partial t} + k(x,u)V(x,u) + LV(x,u)) = 0$$

in the finite time horizon case and

$$\min_u((k(x,u) - \lambda^*)V(x,u) + LV(x,u)) = 0$$

in the infinite time horizon case. One usually needs some technical restrictions on $k$, in particular so that the cost is in fact finite under some control. It has been found more convenient to transform these HJB equations by the logarithmic transformation $\Phi = -\log V$. The $\Phi$ thus defined satisfies the so called Hamilton-Jacobi-Isaacs equation, the counterpart of HJB equation for two person zero sum stochastic differential games, with finite horizon, resp. ergodic payoffs of the type discussed earlier [46], [100]. This transformation, a descendant of the Cole-Hopf transformation that links Burgers equation to the heat equation, was introduced and effectively used by Fleming and his collaborators not only for risk-sensitive control, but also for several interesting spin-offs in large deviations. See, e.g., [50], Ch. 6.

For controlled optimal stopping, the HJB equation gets replaced by the quasi-variational inequalities:

$$\min_u (k(x,u) - c(x,u)V(x) + LV(x,u)) \geq 0,$$
$$h(x) - V(x) \geq 0,$$
$$\min_u (k(x,u) - c(x,u)V(x) + LV(x,u)) (h(x) - V(x)) = 0.$$

These are a slight generalization of variational inequalities appearing in obstacle problems and elsewhere in applied mathematics. The intuition behind these is as follows: If it is optimal not to stop in a neighborhood of $x$, it reduces to the earlier control problem and the HJB equation must hold, i.e., the first inequality above is an equality. If it is optimal to stop at $x$, the minimum cost to go, $V(x)$, must equal the cost on stopping, $h(x)$, i.e., the second inequality above is an equality. In either case, standard dynamic programming heuristic suggests that the appropriate inequality above must hold always. Clearly one of the two equalities must hold at any given point $x$, which leads to the third equality.

The situation for impulse control is similar:

$$\min_u \left(\frac{\partial V}{\partial t}(x,t) + k(x,u) - c(x,u)V(x,t) + LV(x,t,u)\right) \geq 0,$$
$$\min_y(V(y,t) + g(y,x)) - V(x,t) \geq 0,$$
$$\min_u \left(\frac{\partial V}{\partial t}(x,t) + k(x,u) - c(x,u)V(x,t) + LV(x,t,u)\right) \times$$
$$\left(\min_y(V(y,t) + g(y,x)) - V(x,t)\right) = 0.$$

Likewise for optimal switching, we include the control variable '$u$' in the state (thus the value function $V$ has three arguments: $x, u$ and $t$), and consider:

$$\frac{\partial V}{\partial t}(x,u,t) + k(x,u) - c(x,u)V(x,u,t) + LV(x,t,u) \geq 0,$$

$$\min_y(V(x,y,t) + q(y,u)) - V(x,u,t) \geq 0,$$

$$\left(\frac{\partial V}{\partial t}(x,u,t) + k(x,u) - c(x,u)V(x,u,t) + LV(x,u,t)\right) \times$$

$$\left(\min_y(V(x,y,t) + q(y,u)) - V(x,u,t)\right) = 0.$$

See [11], [12] for an extensive treatment of applications of variational and quasi-variational inequalities in stochastic control. A more probabilistic treatment of optimal stopping is found in [106]. See [103] for some recent contributions to optimal switching.

In each case above, the appropriate verification theorem holds. Note also that the verification theorem, coupled with a standard measurable selection theorem (see, e.g., [102]) guarantees an optimal *precise* Markov or stationary Markov control (as applicable). This is because the respective minima are in fact attained at Dirac measures. See [29] for the inequalities for 'stochastic hybrid control'.

### 5.2. HJB equation - the degenerate case

In the degenerate case, the HJB equation typically does not have classical solutions. This has lead to the notion of *viscosity solutions* that provides a unique characterization of the desired solution within a larger class (typically, that of continuous functions). We shall describe this notion in the case of infinite horizon discounted cost problems.

Say that $V$ is a viscosity solution of (24) if for any $\psi \in C^2(\mathcal{R}^d)$,

- at each local maximum of $V - \psi$,

$$\min_u(k(x,u) - c(x,u)V(x) + L\psi(x,u)) \geq 0,$$

and,
- at each local minimum of $V - \psi$,

$$\min_u(k(x,u) - c(x,u)V(x) + L\psi(x,u)) \leq 0.$$

To see why this makes sense in the first place, note that if $V$ were $C^2$, then at each local maximum of $V - \psi$ the gradients of $V, \psi$ would be equal and the Hessian of $V - \psi$ would be negative definite. Thus if $V$ satisfied the HJB equation, $(V, \psi)$ would satisfy the first inequality above at this point. A similar logic applies to the second statement.

Note that if one were to add a term $\nu \Delta V$ to $LV$, $\nu > 0$, in (24), then it would be the HJB equation corresponding to replacing $\sigma(\cdot)$ by $\sqrt{\sigma(\cdot)\sigma^T(\cdot) + \nu I_d}$, $I_d$ being the $d \times d$ identity matrix. This is non-degenerate and thus has a classical solution $V_\nu$ as described above. The viscosity solution is the limit of these as $\nu \downarrow 0$. The term $\nu \Delta V$ appears in equations of fluid mechanics as the 'viscosity' term, hence the terminology. An alternative equivalent definition of viscosity solutions is possible in terms of sub-differentials [89]-[91].

The value function can be shown to be the unique viscosity solution of the HJB equation in an appropriate function class for a wide variety of control problems [50], [89]-[91]. See [55], [92], [113] for the corresponding development in case of variational inequalities.

This leaves open the issue of a verification theorem wherein the utility of this approach finally resides. While this is not as routine as in the non-degenerate case, recent work using non-smooth analysis has made it possible [118].

We mention now two abstractions of the dynamic programming principle which led to the HJB equations. The first is the martingale dynamic programming principle formulated first in [110] (written in 1974, though published much later) and developed further in [34], [104]. For the finite horizon problem above, this reduces to the observation that

$$V(X(t)) + \int_0^t LV(X(s), u(s))ds, \ t \in [0, T],$$

is a submartingale w.r.t. $\{\mathcal{F}_t\}$ and is a martingale if and only if $u(\cdot)$ is optimal. Similar statements can be formulated for the other problems. The second approach is the nonlinear semigroup developed in [101]. This is the semigroup of operators

$$S_t f \stackrel{def}{=} \min E[\int_0^t e^{-\int_0^s c(X(y),u(y))dy} k(X(s), u(s))ds$$
$$+ e^{-\int_0^t c(X(y),u(y))dy} f(X(t))|X(0) = x],$$

where the minimum is over all admissible controls. Under our hypotheses, this can be shown to be a semigroup of positive nonlinear operators $C_b(\mathcal{R}^d) \to C_b(\mathcal{R}^d)$ which is the *lower envelope* of the corresponding Markov (linear) semigroups corresponding to constant controls $u(\cdot) \equiv a \in U$, in a precise sense. The associated infinitesimal generator has the form

$$\hat{L}f = \min_u (Lf(x,u) + k(x,u) - c(x,u)f).$$

The above are resp. the controlled counterparts of the 'martingale problem' and the 'semigroup approach' in Markov process theory, and have the advantage that they generalize naturally to more abstract semimartingale, resp. Markov process control problems.

### 5.3. The stochastic maximum principle

There has also been a considerable body of work on extending the theory of necessary conditions for optimality based on the Pontryagin maximum principle from deterministic optimal control to the stochastic setting. The earliest effort in this direction is perhaps [85]. It may be recalled that the maximum principle involves an 'adjoint process' which evolves backward in time with a given terminal condition. Since stochastic control comes with the additional baggage of the 'arrow of time' specified by the increasing filtration and associated adaptedness / nonanticipativity issues, this extension is nontrivial and much hard work went into it. See, e.g., [64], which was a landmark contribution in this domain, and the references therein. The advent of 'backward stochastic differential equations' provided a natural framework for handling this, culminating in the very general stochastic maximum principle (for the finite horizon problem) reported in [116]. A typical b.s.d.e. is of the form

$$dY(t) = h(t, Y(t), Z(t))dt + Z(t)dW(t), \ t \in [0, T],$$

with the terminal condition $Y(T) = \xi$. Here, for $\mathcal{F}_t^W \stackrel{def}{=}$ the natural filtration of $W(\cdot)$, $\xi$ is a prescribed square integrable random variable measurable with respect to $\mathcal{F}_T^W$. A solution is a pair of $\{\mathcal{F}_t^W\}$–adapted processes $Y(\cdot), Z(\cdot)$ such that

$$Y(t) = \xi - \int_t^T h(t, Y(s), Z(s))ds - \int_t^T Z(s)dW(s), \ t \in [0, T].$$

Under a Lipschitz condition on $h$, a unique solution can be shown to exist in a suitable class of $\{\mathcal{F}_t^W\}$–adapted processes ([116], Chapter 7). See [95] for an extensive account of coupled forward-backward stochastic differential equations and their applications to stochastic control and mathematical finance. See also [42].

A special case of the stochastic maximum principle, for $\sigma$ independent of control, is as follows. Assume that $m, \sigma, k, h$ are bounded, twice continuously differentiable in the space $(x)$ variable with the first and second order partial derivatives satisfying the Lipschitz condition. We confine ourselves to $\{\mathcal{F}_t^W\}$–adapted controls $u(\cdot)$. Let $p(\cdot), q(\cdot) = [q_1(\cdot)|q_2(\cdot)|\cdots|q_d(\cdot)]$ be processes adapted to the natural filtration of $W(\cdot)$ and satisfying the backward stochastic differential equation

$$\begin{aligned} dp(t) &= -(\nabla_x m(X(t), u(t))^T p(t) + \sum_i \nabla_x \sigma^i(X(t))^T q_i(t) \\ &\quad - \nabla_x k(X(t), u(t)))dt + q(t)dW(t), \ t \in [0, T], \end{aligned} \quad (25)$$

with the terminal condition $p(T) = -\nabla_x h(x(T))$. Here $\sigma^j(\cdot)$ denotes the $j$-th column of $\sigma(\cdot)$. Under stated conditions, (25) can be shown to have an a.s. unique solution $(p(\cdot), q(\cdot))$. The process $p(\cdot)$ is the desired adjoint process. The maximum principle then states that if $(X(\cdot), u(\cdot))$ is an optimal pair, then one must have

$$\langle p(t), m(X(t), u(t))\rangle - k(X(t), u(t)) = \max_u \left(\langle p(t), m(X(t), u)\rangle - k(X(t), u)\right) \tag{26}$$

for a. e. $t \in [0, T]$. In fact, the full statement of the stochastic maximum principle in [116] is much more general, allowing for a controlled diffusion matrix $\sigma$.

Comparing with the verification theorem of dynamic programming, one would expect $p(t)$ above to correspond to $-\nabla_x V(X(t), t)$. This may be shown under very strong conditions. More generally, a relationship along these lines has been established in [118] (see also [116]).

### 5.4. Partial observations

The dynamic programming principle under partial observations is usually derived by moving over to the 'separated' control problem of controlling the associated nonlinear filter. In the simpler cases, the 'integral' form of the dynamic programming principle is easy to justify. For example, for the finite horizon problem, define the value function

$$V(\pi, t) \stackrel{def}{=} \min E[\int_t^T \pi_s(k(\cdot, u(s)))ds + \pi_T(h)|\pi_t = \pi],$$

where the minimum is over all wide sense admissible controls. This satisfies: for $\Delta > 0$,

$$V(\pi_t, t) = \min E[\int_t^{t+\Delta} \pi_s(k(\cdot, u(s)))ds + V(\pi_{t+\Delta}, t+\Delta)|\pi_t = \pi],$$

with the minimum attained if and only if $u(s)|_{s\in[t,t+\Delta]}$ is an optimal choice. Analogous statements can be made for the unnormalized law as the state variable. To get a 'differential' form of this principle in terms of an HJB equation is hard, as the state space, $\mathcal{P}(\mathcal{R}^d)$ or $\mathcal{M}(\mathcal{R}^d)$ (alternatively, the more popular $L_2(\mathcal{R}^d)$ when a square integrable density for the conditional law is available), is infinite dimensional. This has been approached through the theory of viscosity solutions for infinite dimensional p.d.e.s [57], [90]. As for the more abstract versions, the martingale formulation of the dynamic programming principle for the separated control problem is a straightforward counterpart of the completely observed case. See, however, [34] for a different development which derives a martingale dynamic programming principle in a very general set-up (see also [52]). The Nisio semigroup has been developed in [23], [45]. See [32], [112] for some recent developments in the stochastic maximum principle under partial observations and [7] for an early contribution.

### 6. Computational issues

Stochastic control problems with elegant explicit solutions tend to be few. There are, however, some notable exceptions of great practical importance, such as the celebrated 'Linear-Quadratic-Gaussian' problem with linear state dynamics and quadratic cost, which has become standard textbook material [31]. More often than not the controlled diffusion problems do not lead to explicit analytic solutions and one has to resort to approximations and numerical computations. This has led to much research in approximation and computational issues. We briefly survey some of the main strands of this research.

One popular method has been to consider controlled Markov chain approximations to controlled diffusions, thereby moving over to discrete time and discrete state space. One then analyzes the resulting discrete problem by standard schemes available for the same. See [84] for an extensive account of a rigorous theory for this well developed approach. [83] contains some recent extensions of this approach to stochastic differential games.

The other important approach considers the infinite linear program implicit in the occupation measure based approach and uses linear programming tools (see, e.g., [33], [66]). The ensuing linear program, however, is infinite dimensional and its approximation by a finite linear program is needed [97].

The HJB equation, being a nonlinear p.d.e., is open to numerical techniques developed for the same. The most important recent developments on this front are the ones propelled by the viscosity solutions revolution that use stability results for viscosity solutions for rigorous justification. See, e.g., [4].

The recent developments in simulation-based approximate dynamic programming [14], however, have not caught on in controlled diffusion literature to a large extent, but there is considerable interest in the finance community for such 'Monte Carlo' techniques - see, e.g., [56].

For numerical analysis of stochastic differential equations in general, [72] is the standard source. A good source for 'Monte Carlo' for diffusions is [87].

## 7. Other problems

Here we list some other subareas of controlled diffusions that will not be discussed at length here. Only a brief description is given, with some representative references.

1. *Singular control:* These are problems involving an additive control term in the stochastic 'integral' equation that is of bounded variation, but not necessarily absolutely continuous with respect to the Lebesgue measure. That is,

   $$X(t) = X_0 + \int_0^t m(X(s))ds + A(t) + \int_0^t \sigma(X(s))dW(s), \ t \geq 0,$$

   where $A(\cdot)$ is the control. Typically it can be 'local time at a boundary' that confines the process to a certain bounded region. This research originated in heavy traffic limits of controlled queues [60], [107]. See [2], [35], [51], [74], [75] for some recent contributions and applications to finance.

2. *Adaptive control:* This concerns the situation when the exact model of the controlled system is not known and has to be 'learnt on line' while controlling it. Several alternative approaches to this problem exist in the discrete time stochastic control literature, but the only one that seems to have been followed to any significant extent in controlled diffusions is the 'self-tuning' control [22], [36]. In this, one enforces a separation of estimation and control by estimating the model by some standard statistical scheme (usually parametric), and at each time using the control choice that would be optimal for that time and state if the current estimate were indeed the correct model. This runs into the usual 'identifiability' problem: several models may lead to control choices that in turn lead to identical output behavior, making it impossible to discriminate between these models. Many variations have been suggested to work around this problem, such as additional randomization of controls as 'probes'.
3. *Control of modified diffusions and control with additional constraints:* Issues similar to those of the preceding section have been explored for reflected diffusions [30] (which often arise as heavy traffic approximation of controlled queues [61]), diffusions with 'jumps' or switching modes [54], [88], etc. Another related development is control under additional constraints [27]. Here the controller seeks to minimize one cost functional subject to a bound on one or more ancillary cost functionals.
4. *Multiple timescales:* These are problems wherein different components of the controlled diffusion move on different time-scales, as in:

$$dX_1(t) = m^{(1)}(X_1(t), X_2(t), u(t))dt + \sigma^{(1)}(X_1(t), X_2(t))dW_1(t),$$
$$dX_2(t) = \frac{1}{\epsilon}m^{(2)}(X_1(t), X_2(t), u(t))dt + \frac{1}{\sqrt{\epsilon}}\sigma^{(2)}(X_1(t), X_2(t))dW_2(t),$$

where $\epsilon > 0$ is 'small'. This implies in particular that $X_2(\cdot)$ operates on a much faster time-scale than $X_1(\cdot)$. Intuitively, one would expect $X_2(\cdot)$ to see $X_1(\cdot)$ as quasi-static, whereas $X_1(\cdot)$ sees $X_2(\cdot)$ as almost equilibrated. This intuition is confirmed by the analysis which allows one to analyze $X_1(\cdot)$ with its dynamics averaged over the asymptotic behavior (read 'stationary distribution' in the asymptotically stationary case) of $X_2(\cdot)$ *when* the latter is analyzed by freezing the $X_1(\cdot)$ in *its* dynamics as though it were a constant parameter [69], [82].
5. *Game problems:* These are the problems that involve more than one controller with possibly different costs. The simplest is the two person zero sum case where two controllers have cost functionals that sum to zero, i.e., the cost of one is the reward of the other. The key result in this case is the minmax theorem which establishes the existence of a value, equalling both the minimum of the maximum (over the opponent's choices) cost paid by the first and the maximum of the minimum (over the opponent's choices) reward gained by the other. This then is characterized by the appropriate Hamilton-Jacobi-*Isaacs* equation for the value function, which corresponds to replacing the 'min' in the HJB equation by 'minmax' or 'maxmin'. The more general $N$−person noncooperative case has $N$ controllers with different cost functionals. This is more complicated and one looks for a *Nash equilibrium*, i.e., a control policy profile for the controllers whereby no single controller can lower her cost by choosing differently if the rest don't change their controls. This leads to a coupled system of HJB equations, coupled through the minimizing controls of each other. See [10], [28], [114] for a sampler.

6. *Mathematical finance:* This has proved to be a rich source of problems in stochastic control in recent years, e.g., in option pricing, portfolio optimization, etc. We have already seen some examples in Section 3. The area is still exploding and merits a separate full length review. See [70], [71], [105] for a perspective and [49], [58], [68] for a sample of recent contributions.

What next? To mention a few of the current themes, the most prominent of course remain the problems emerging from mathematical finance and heavy traffic limits of queues. Risk-sensitive control is another area which still offers interesting open problems, as are control of degenerate diffusions and control under partial observations. Extensions to infinite dimensional problems also present several challenges of a technical nature. The biggest challenge, however, is on the computational front. Fast and accurate computational schemes are sought in particular by the finance community.

<